\documentclass{article}
\usepackage{listings}
\usepackage[T1]{fontenc}

\usepackage[english]{babel}
\usepackage{amssymb}
\usepackage{amsthm}

\newtheorem{theorem}{Theorem}

\newcommand{\Iso}{\ensuremath{\mathit{Iso}}}

\newcommand \n {\{ 1, 2, \ldots , n \}}

\lstdefinelanguage{GAP}{%
	morekeywords={%
		Assert,Info,IsBound,QUIT,%
		TryNextMethod,Unbind,and,break,%
		continue,do,elif,%
		else,end,false,fi,for,%
		function,if,in,local,%
		mod,not,od,or,%
		quit,rec,repeat,return,%
		then,true,until,while%
	},%
	sensitive,%
	morecomment=[l]\#,%
	morestring=[b]",%
	morestring=[b]',%
}[keywords,comments,strings]

\makeatletter
\def\@fnsymbol#1{\ensuremath{\ifcase#1\or *\or **\or \ddagger\or
		\mathsection\or \mathparagraph\or \|\or **\or \dagger\dagger
		\or \ddagger\ddagger \else\@ctrerr\fi}}
\makeatother
\begin{document}

\title{Improvements for lower bounds of mutually orthogonal Latin squares of sizes $54$, $96$ and $108$}

\author{
  R. Julian R. Abel \thanks{r.j.abel@unsw.edu.au } \\
School of Mathematics and Statistics\\
UNSW Sydney\\ NSW 2052, Australia\\
\and	
Ingo Janiszczak\thanks{ingo.janiszczak@uni-due.de} \ \ and Reiner Staszewski\\
Faculty of Mathematics\\
University of Duisburg-Essen\\
45127 Essen, Germany\\
 }

\maketitle

\begin{abstract}
{
We will show that there are at least 8, 10 and 9 mutually orthogonal Latin squares (MOLS) of orders $n=54$, $96$ and $108$.  The cases 
$n=54$ and $96$ are obtained by constructing separable permutation codes
consisting of $8 \times 54$ and $10 \times 96$ codeword respectively; 
in addition, these  codes  respectively  have  lengths $54$, $96$
and minimum distances $53$,  $95$. 
Here we will follow exactly the procedure given in \cite{JS2019}.  The 
case $n=108$ is obtained by constructing a $(108,10,1)$  difference matrix. Also, an error in \cite{ACD} for $n=45$  will be corrected.
}
\end{abstract}



\section{Introduction}
\label{sect:intro}

Let $n$ be a positive integer and let $V$ and $W$ be sets consisting of $n$ elements.
A Graeco-Latin square of order $n$  is a $n \times n$ array 
$M = ((v_{ij},w_{ij}))$ with entries in 
the Cartesian product $V \times W$ such that the set of all different 
entries in $M$ is $V \times W$ and the arrays $L_1 = (v_{ij})$ and
$L_2 = (w_{ij})$ are Latin squares, i.e. all rows and all columns of $L_1$
and $L_2$ are permutations of $V$ and $W$ respectively. 
Two such Latin squares are called orthogonal if the array $((v_{ij},w_{ij}))$
is a Graeco-Latin square. A set of $t$ mutually orthogonal Latin squares (MOLS) 
of order $n$ (or $t$ MOLS$(n)$) consists of $t$ Latin squares of order $n$ 
which are pairwise orthogonal.  Let $N(n)$ be  size of the largest set of MOLS of order $n$.

In this paper, we obtain improvements for N$(n)$ when $n \in \{54,96,108\}$. 
These MOLS determined using two different methods. For $n=54$ and $n=96$, permutation codes 
corresponding to the relevant MOLS were determined. In the case of $n=108$, a difference 
matrix with $10$ rows was constructed; from this difference matrix,  $9$ MOLS for $n=108$ can be obtained. 

Since parts of the present paper for the cases $n=54$ and $n=96$ are to be seen as a supplement 
to \cite{JS2019}, we will simply extract the relevant part for the necessary notations and 
definitions.

Let $X$ be an arbitrary finite set of order $n$.
Let $S_X := \{ \sigma : X \rightarrow X \mid \sigma \;\mbox{is bijective} \}.$
This group acts on $X$ from the right, and
for $x \in X, \sigma \in S_X,$ we denote the image of $x$ under $\sigma$ by $x^{\sigma}.$
If $X = \n$, then  we will write $S_n$ instead of $S_X.$
For $\sigma, \tau \in S_X$ the Hamming distance $d_H(\sigma, \tau) = d(\sigma, \tau)$ equals $n$ minus 
the number of fixed points of $\sigma^{-1} \tau$.
A subset $C$ of $S_X$ is called a \emph{permutation code} 
or \emph{permutation array} of length $n$ and of minimum distance $d = d(C)$ if
$$d(C) := \min \ \{d_H(\sigma, \tau) \mid \sigma, \tau \in C , \sigma \ne \tau \} \ .$$

Such a permutation array $C$ is called an $(n,d)$-PA.
If $C$ is an $(n,n-1)$-PA, then  $C$ is called $(r,m)$-separable if it is the join of $m$ disjoint $(n,n)$-PAs 
$P_1, \ldots, P_m$ of cardinality $r$  such that 
for all pairs $P_s, P_t$, $s \ne t$, the distance $d_H(\sigma, \tau)$ equals $n-1$ for 
all $\sigma \in P_s$ and all $\tau \in P_t$.   When $r = n$, existence of an  $(n,m)$-separable 
$(n,n-1)$-PA is equivalent to existence of $m$ MOLS of order $n$ \cite{CKL04}.
In this case we order the elements in $P_s$ such that $P_s = \{ \sigma_{sk} : 1 \le k \le n\}$ for all 
$1 \le s \le m$. From this $m$ MOLS $A_1, A_2, \ldots, A_m$ of order $n$ can be obtained as follows:
$A_{s}(i,j) = k$ if $j^{\sigma_{sk}} = i$.

$S_n$ is a metric space via the Hamming distance, and the isometry group $\Iso(n)$ 
is isomorphic to the wreath product $S_n \wr S_2$ \cite{F60}.
$\Iso(n)$ can be described as a subgroup of $S_{2n}$.
Let $B_1$ and $B_2$ be the naturally embedded subgroups isomorphic to $S_n$
acting on the sets $\{1,2,\ldots ,n\}$ and $\{n+1,n+2,\ldots ,2n\}$, respectively. 
Moreover let $t_n := (1,n+1)(2,n+2)\cdots (n,2n) \in S_{2n}$. Then 
$\Iso(n) = \langle B_1, B_2, t_n\rangle = (B_1 \times B_2):\langle t_n \rangle$,
and the action of this group on $B_1$ from right is given by
$$b * x := \left\{ \begin{array}{lll}
x^{-1} b \; & \mbox{if} \; & x \in B_1 \\
b \varphi(x) \; & \mbox{if} \; & x \in B_2 \\
b^{-1} \; & \mbox{if} \; & x = t_n, \\
\end{array}
\right. $$
where $b \in B_1$, $x \in Iso(n)$ and $\varphi$ denotes the natural isomorphism from $B_2$ to $B_1$.
Moreover, if $U$ is a subgroup of $\Iso(n)$ and $b \in B_1,$  then
$b * U := \{ b * u \mid u \in U \}$ denotes the $U-$orbit of $b.$
Our codes $C$ will now be regarded as subsets of $B_1$, and a code $C$ will be called 
$U-$invariant if $C$ is closed under the action of $U$. This means that $C$ is a union of $U-$orbits.
For the strategy of constructing $(n,m)$-separable PAs that are invariant under a given subgroup $U$ of
$\Iso(n)$ we refer to \cite{JLOS2015}. But now we will only join separable $U$-orbits.

In section~\ref{sect:impr}  we will prove the existence of $8$ MOLS for $n=54$ and $10$ MOLS for  
$n=96$ by constructing a $(54,8)$-separable PA and a $(96,10)$-separable PA.
These are invariant under specific subgroups  $U$ of $Iso(n)$ using 
the described group action of $Iso(n)$ on $S_n$ and will be
given by generators of $U$ and representatives of the $U$-orbits. 

As all the data for the constructed PAs is contained in the symmetric group $S_{2n}$, 
it can easily be read into a computer algebra system such as MAGMA \cite{BCP97} or GAP \cite{GAP}.
The orbits can then be computed using elementary programs consisting of a few lines of 
GAP or MAGMA-code.
In the appendix we give a GAP code for calculating the orbits, the resulting PA and the MOLS.
There is also a routine to check if the PA is $(n,m)$-separable.

Now we  define a difference matrix.
Let $G$ be a finite  group of order $n$ and let $k\geq 2$ be an integer. An $(n,k,1)$-difference matrix (DM) over $G$, or briefly $(n,k,1)$-DM over $G$, is a $k\times n$
matrix $D=(d_{i,j})$ with entries from $G$ such that  for any two distinct rows $s$ and $t$, the multiset of differences $\{(d_{s,j})(d_{t,j})^{-1}:1\leq j\leq n\}$  contains each element of $G$ exactly once.  For most known MOLS obtained from difference matrices, the group $G$ is abelian and the expression $d_{s,j} - d_{t,j}$
instead of $(d_{s,j})(d_{t,j})^{-1}$  is  more commonly used. 

If $D$ is an $(n,k,1)$-DM over $G$, and the first row of $D$ is normalised (that is, 
for all $j$,  the entries in column $j$ have been multiplied on the right 
by $d_{1,j}^{-1}$  so that the first entry in each column is $1$, the identity 
element in $G$) then $k-1$ MOLS$(n)$,  $L_1, L_2, \ldots, L_{k-1}$ can be obtained 
as follows.  First let $G= \{g_1 = 1, g_2, \ldots, g_n\}$.  Then for $j=1,2, \ldots, n$,
the $j$'th row of $L_m$ is obtained by multiplying the $(m+1)'$th row of $D$ on the right
by $g_j$.   In section~\ref{sect:impr} we obtain a $(108,10,1)$-DM which gives $9$ MOLS$(108)$.

\section{Improvements}
\label{sect:impr}

\bigskip

Let $n = 54$ and $U := \langle x_1,x_2,x_3,x_4,x_5 \rangle \le Iso(54)$ of order $243$ 
generated by $x_1,x_2,x_3,x_4,x_5$ given in cycle structure

\noindent
$
x_1:=(1, 6, 26)(2, 4, 27)(3, 5, 25)(7, 10, 13)(8, 11, 14)(9, 12, 15)(16, 21, 23)(17, 19, 24)(18, 20, 22)(28, 33, 53)\\
 \ (29, 31, 54)(30, 32, 52)(34, 37, 40)(35, 38, 41)(36, 39, 42)(43, 48, 50)(44, 46, 51)(45, 47, 49)(55, 60, 80)\\
 \ (56, 58, 81)(57, 59, 79)(61, 64, 67)(62, 65, 68)(63, 66, 69)(70, 75, 77)(71, 73, 78)(72, 74, 76)(82, 87, 107)\\
 \ (83, 85, 108)(84, 86, 106)(88, 91, 94)(89, 92, 95)(90, 93, 96)(97, 102, 104)(98, 100, 105)(99, 101, 103),\\
x_2:= (1, 11, 19)(2, 12, 20)(3, 10, 21)(4, 15, 22)(5, 13, 23)(6, 14, 24)(7, 16, 25)(8, 17, 26)(9, 18, 27)\\
 \ (28, 38, 46)(29, 39, 47)(30, 37, 48)(31, 42, 49)(32, 40, 50)(33, 41, 51)(34, 43, 52)(35, 44, 53)(36, 45, 54),\\
x_3:=(1, 2, 3)(4, 5, 6)(7, 8, 9)(10, 11, 12)(13, 14, 15)(16, 17, 18)(19, 20, 21)(22, 23, 24)(25, 26, 27)\\
 \ (28, 29, 30)(31, 32, 33)(34, 35, 36)(37, 38, 39)(40, 41, 42)(43, 44, 45)(46, 47, 48)(49, 50, 51)(52, 53, 54),\\
x_4:=(55, 65, 73)(56, 66, 74)(57, 64, 75)(58, 69, 76)(59, 67, 77)(60, 68, 78)(61, 70, 79)(62, 71, 80)(63, 72, 81)\\
 \ (82, 92, 100)(83, 93, 101)(84, 91, 102)(85, 96, 103)(86, 94, 104)(87, 95, 105)(88, 97, 106)(89, 98, 107)\\
 \ (90, 99, 108),\\
x_5:=(55, 56, 57)(58, 59, 60)(61, 62, 63)(64, 65, 66)(67, 68, 69)(70, 71, 72)(73, 74, 75)(76, 77, 78)(79, 80, 81)\\
 \ (82, 83, 84)(85, 86, 87)(88, 89, 90)(91, 92, 93)(94, 95, 96)(97, 98, 99)(100, 101, 102)(103, 104, 105)\\
 \ (106, 107, 108).
$

Furthermore let

\noindent
$
a_1:=(1, 34, 20, 53, 10, 45)(2, 51, 12, 22, 43, 23, 9, 47, 26)(3, 13, 35, 14, 25, 37, 18, 21, 42)\\
 \ (4, 17, 28, 7, 11, 32, 19, 6, 54)(5, 46, 15, 30, 24, 39)(8, 41)(16, 50)(27, 31)(29, 48, 38)(36, 44, 52),\\
a_2:=(1, 4, 21, 15, 27, 5, 30, 19, 16, 2, 7, 23, 18, 50, 29, 28, 36, 49, 14, 48, 45, 17, 13, 10, 25, 40, 33, 47)\\
 \ (3, 43, 42, 39, 54, 44, 53, 31, 41, 24)(6, 9, 12, 22, 37, 46, 38, 11, 34, 35, 26, 20, 52, 8, 32),\\
a_3:=(1, 38, 26, 35, 6, 41)(2, 39, 27, 36, 4, 42)(3, 37, 25, 34, 5, 40)(16, 30, 51, 18, 29, 50, 17, 28, 49)\\
 \ (19, 33, 45, 21, 32, 44, 20, 31, 43)(22, 54, 48, 24, 53, 47, 23, 52, 46),\\
a_4:=(1, 41, 37, 19, 43, 8, 28, 9, 48, 35, 4, 22, 34, 13, 38, 10)(2, 32, 6, 5, 14, 30, 16, 17, 7, 39)\\
 \ (3, 49, 11, 21, 53, 50)(15, 47, 45, 18, 27, 31)(20, 36, 24, 44, 25, 51)(23, 54, 42, 29, 26, 40, 46, 52, 33),\\
a_5:=(1, 25, 11, 23, 8, 34, 35, 43, 30, 32, 41, 36, 54)(2, 9, 45, 46, 48, 38, 22, 27)\\
 \ (3, 17, 5, 33, 50, 7, 53, 20, 47, 29, 51, 18, 15, 12, 4, 49, 26, 21, 28, 42, 44, 39, 6, 40, 52, 10, 14)(13, 19, 37)(16, 24),\\
a_6:=(1, 54, 48, 40, 20, 34, 6, 38, 51, 41, 52, 36, 53, 15, 4, 19, 12, 21, 16, 14, 10, 44, 43, 28, 17, 25, 18, 33, 3)\\
 \ (2, 9, 8, 22, 47, 27, 5, 13, 30, 46, 29, 31, 35, 39, 37, 7, 42, 49, 11, 26, 50, 45, 23),\\
a_7:=(1, 23, 5, 14, 52, 36, 17, 8, 25, 12, 11, 6, 45, 27, 30, 19, 15, 9, 38, 28, 48, 10, 49, 50, 44, 29, 51,\\
 \ 24, 35, 46, 39, 31, 33, 4, 26, 18, 47, 40, 34, 43, 53, 37, 2, 3, 41, 13, 22, 16)(7, 21, 32, 54),\\
a_8:=(1, 35, 17, 45, 25, 11, 24, 51, 12, 43, 47)(2, 29, 16, 33, 52, 9, 49, 4, 32, 14, 18)(3, 15, 46, 8, 21, 54, 22, 26, 40)\\
 \ (5, 7, 34, 6, 39, 41, 23, 30, 10, 38, 20)(13, 42, 36, 37, 27, 53, 28, 31, 19, 48, 50, 44),\\
a_9:=(1, 27, 29, 21, 26, 36, 8, 23, 44, 54, 52, 14, 32, 5, 31)(2, 47, 39, 18, 24, 11, 40)(3, 42, 51)\\
 \ (4, 12, 25, 22, 33, 50, 49)(6, 45, 13, 10, 48, 16, 35, 15, 43, 9, 34, 53, 7, 30, 37, 19, 46, 20, 41)(17, 28),\\
a_{10}:=(1, 18, 34, 33, 12, 7, 36, 51)(2, 27, 44, 49, 11, 25, 53, 32)(3, 9, 54, 40, 10, 16, 43, 42)\\
 \ (4, 28, 14, 29, 24, 37, 15, 39)(5, 38, 22, 48, 23, 30, 6, 46)(8, 45, 31, 20, 17, 52, 50, 19)(13, 47)(21, 26, 35, 41),\\
a_{11}:=(1, 32, 13, 53, 10, 24, 40, 22, 18, 29, 52, 41, 28, 19, 51, 37, 14, 33, 47, 44)\\
 \ (2, 30, 4, 26, 21, 5, 45, 16, 54, 27, 38, 3, 15, 7, 43, 31, 17, 11, 42, 25, 35, 20, 46, 12, 39, 36, 8)\\
 \ (6, 49, 9, 48, 23, 34, 50),\\
a_{12}:=(1, 17, 38, 45, 12, 25, 37, 54)(2, 8, 47, 53, 11, 7, 29, 43)(3, 26, 30, 34, 10, 16, 46, 35)\\
 \ (4, 24, 15, 13, 5, 14, 23, 22)(9, 39, 36, 21)(18, 28, 52, 19, 27, 48, 44, 20)(31, 32, 50, 33, 40, 42, 49, 41),\\
a_{13}:=(1, 15, 51, 31, 46, 18, 34, 53, 44, 50, 38, 37, 45, 35, 5, 41, 28, 36, 40, 19, 20, 4, 17, 21, 47, 54,\\
 \ 23, 8, 25, 3, 2, 29, 30, 27, 16, 7, 52, 32, 42, 49)(6, 14, 26, 43, 13, 24, 33, 11, 39, 9, 10, 22),\\
a_{14}:=(1, 38, 21, 11, 7, 22, 6, 32, 46, 49, 34, 27, 44, 30, 54, 41, 15)(2, 37, 19, 10, 8, 24, 4, 31, 47, 51,\\
 \ 35, 26, 45, 29, 52, 40, 13, 3, 39, 20, 12, 9, 23, 5, 33, 48, 50, 36, 25, 43, 28, 53, 42, 14)(16, 18),\\
a_{15}:=(1, 21, 19, 43, 2, 45, 49, 36, 3, 23, 30, 14, 4, 25, 42, 37, 52, 51, 18, 40, 16, 11, 22, 6, 29, 53)\\
 \ (7, 9, 12, 20, 24, 26, 10, 44, 47, 15, 28, 33, 38, 13, 27)(8, 41, 34, 48, 31, 17)(32, 35, 50, 39)(46, 54),\\
a_{16}:=(1, 42, 5, 44, 18, 48)(2, 34, 28, 6, 47, 50, 25, 33, 36)(3, 17, 4)(7, 41, 45, 12, 43, 38, 14, 29, 32)\\
 \ (8, 22, 21)(9, 37, 19, 31, 23, 35)(10, 26, 15)(11, 49, 13, 53, 27, 30)(16, 51, 54, 20, 52, 46, 24, 39, 40).
$

For $i \in \{1, \cdots ,16\}$ the $a_i$ are codewords of $B_1$ and the orbits $a_i * U$ all have cardinality $27$.
The $(54,53)$ - PA $C := \bigcup_{i=1}^{16} a_i * U $ is $(54, 8)$ - separable.
This proves
\bigskip
\begin{theorem}
\label{thm:8MOLS54}
$N(54) \ge 8$.
\end{theorem}

\bigskip 

Let $n = 96$.
In \cite{JS2019} the existence of 8 MOLS has been proven by constructing
a $(96, 8)$ - separable code. This implies $N(96) \ge 8$. Now we like to improve this bound using
a different isometry group.
Let $U := \langle x_1,x_2,x_3,x_4 \rangle \le Iso(96)$ of 
order $2304$ generated by $x_1,x_2,x_3,x_4$ given in cycle structure

\noindent
$
x_1:=(1, 4)(2, 3)(5, 16)(6, 15)(7, 13)(8, 14)(9, 11)(10, 12)(17, 20)(18, 19)(21, 32)(22, 31)\\
 \ (23, 29)(24, 30)(25, 27)(26, 28)(33, 36)(34, 35)(37, 48)(38, 47)(39, 45)(40, 46)(41, 43)\\
 \ (42, 44)(49, 52)(50, 51)(53, 64)(54, 63)(55, 61)(56, 62)(57, 59)(58, 60)(65, 68)(66, 67)\\
 \ (69, 80)(70, 79)(71, 77)(72, 78)(73, 75)(74, 76)(81, 84)(82, 83)(85, 96)(86, 95)(87, 93)\\
 \ (88, 94)(89, 91)(90, 92)(97, 99)(98, 100)(101, 111)(102, 112)(103, 110)(104, 109)(105, 108)\\
 \ (106, 107)(113, 115)(114, 116)(117, 127)(118, 128)(119, 126)(120, 125)(121, 124)(122, 123)\\
 \ (129, 131)(130, 132)(133, 143)(134, 144)(135, 142)(136, 141)(137, 140)(138, 139)(145, 147)\\
 \ (146, 148)(149, 159)(150, 160)(151, 158)(152, 157)(153, 156)(154, 155)(161, 163)(162, 164)\\
 \ (165, 175)(166, 176)(167, 174)(168, 173)(169, 172)(170, 171)(177, 179)(178, 180)(181, 191)\\
 \ (182, 192)(183, 190)(184, 189)(185, 188)(186, 187),\\
x_2:=(1, 20, 33, 4, 17, 36)(2, 19, 34, 3, 18, 35)(5, 32, 37, 16, 21, 48)(6, 31, 38, 15, 22, 47)\\
 \ (7, 29, 39, 13, 23, 45)(8, 30, 40, 14, 24, 46)(9, 27, 41, 11, 25, 43)(10, 28, 42, 12, 26, 44)\\
 \ (49, 68, 81, 52, 65, 84)(50, 67, 82, 51, 66, 83)(53, 80, 85, 64, 69, 96)(54, 79, 86, 63, 70, 95)\\
 \ (55, 77, 87, 61, 71, 93)(56, 78, 88, 62, 72, 94)(57, 75, 89, 59, 73, 91)(58, 76, 90, 60, 74, 92)\\
 \ (97, 119, 129, 103, 113, 135)(98, 120, 130, 104, 114, 136)(99, 126, 131, 110, 115, 142)\\
 \ (100, 125, 132, 109, 116, 141)(101, 123, 133, 107, 117, 139)(102, 124, 134, 108, 118, 140)\\
 \ (105, 128, 137, 112, 121, 144)(106, 127, 138, 111, 122, 143)(145, 167, 177, 151, 161, 183)\\
 \ (146, 168, 178, 152, 162, 184)(147, 174, 179, 158, 163, 190)(148, 173, 180, 157, 164, 189)\\
 \ (149, 171, 181, 155, 165, 187)(150, 172, 182, 156, 166, 188)(153, 176, 185, 160, 169, 192)\\
 \ (154, 175, 186, 159, 170, 191),\\
x_3:=(1, 96, 17, 64, 33, 80)(2, 87, 29, 63, 41, 75)(3, 86, 19, 54, 35, 70)(4, 94, 24, 53, 44, 74)\\
 \ (5, 92, 26, 52, 46, 72)(6, 83, 22, 51, 38, 67)(7, 93, 31, 57, 43, 66)(8, 85, 28, 58, 36, 78)\\
 \ (9, 91, 18, 55, 45, 79)(10, 84, 30, 56, 37, 76)(11, 82, 23, 61, 47, 73)(12, 90, 20, 62, 40, 69)\\
 \ (13, 95, 25, 59, 34, 71)(14, 88, 21, 60, 42, 68)(15, 89, 27, 50, 39, 77)(16, 81, 32, 49, 48, 65)\\
 \ (97, 164, 137, 151, 125, 192)(98, 168, 130, 152, 114, 184)(99, 162, 140, 158, 120, 182)\\
 \ (100, 173, 132, 157, 116, 189)(101, 171, 133, 155, 117, 187)(102, 175, 142, 156, 122, 179)\\
 \ (103, 169, 139, 145, 128, 181)(104, 166, 131, 146, 124, 190)(105, 167, 141, 160, 113, 180)\\
 \ (106, 163, 134, 159, 126, 188)(107, 161, 144, 149, 119, 185)(108, 174, 136, 150, 115, 178)\\
 \ (109, 176, 129, 148, 121, 183)(110, 172, 138, 147, 118, 191)(111, 170, 143, 154, 127, 186)\\
 \ (112, 165, 135, 153, 123, 177),\\
x_4:=(1, 86, 17, 54, 33, 70)(2, 94, 29, 53, 41, 74)(3, 96, 19, 64, 35, 80)(4, 87, 24, 63, 44, 75)\\
 \ (5, 89, 26, 50, 46, 77)(6, 81, 22, 49, 38, 65)(7, 88, 31, 60, 43, 68)(8, 95, 28, 59, 36, 71)\\
 \ (9, 90, 18, 62, 45, 69)(10, 82, 30, 61, 37, 73)(11, 84, 23, 56, 47, 76)(12, 91, 20, 55, 40, 79)\\
 \ (13, 85, 25, 58, 34, 78)(14, 93, 21, 57, 42, 66)(15, 92, 27, 52, 39, 72)(16, 83, 32, 51, 48, 67)\\
 \ (97, 165, 137, 153, 125, 177)(98, 170, 130, 154, 114, 186)(99, 175, 140, 156, 120, 179)\\
 \ (100, 171, 132, 155, 116, 187)(101, 173, 133, 157, 117, 189)(102, 162, 142, 158, 122, 182)\\
 \ (103, 167, 139, 160, 128, 180)(104, 163, 131, 159, 124, 188)(105, 169, 141, 145, 113, 181)\\
 \ (106, 166, 134, 146, 126, 190)(107, 176, 144, 148, 119, 183)(108, 172, 136, 147, 115, 191)\\
 \ (109, 161, 129, 149, 121, 185)(110, 174, 138, 150, 118, 178)(111, 168, 143, 152, 127, 184)\\
 \ (112, 164, 135, 151, 123, 192).
$

Furthermore let

\noindent
$
a_1:=(2, 9, 30, 43, 34, 77, 51, 76, 54, 70, 17, 3, 58, 15, 16, 7, 19, 87, 72, 39, 61, 49, 44, 42, 59, 55)\\
 \ (4, 50, 35, 31, 6, 63, 46, 47, 68, 32, 14, 83, 5, 56, 10, 22, 81, 88, 79, 85)\\
 \ (8, 28, 37, 18, 12, 94, 29, 36, 24, 71, 48, 75, 62, 57)(11, 86)\\
 \ (13, 92, 20, 96, 82, 95, 90, 74, 33, 69, 25, 73, 41, 52, 67, 23, 80, 93, 21, 89, 66, 91, 27, 45,\\
 \ 40, 53, 78, 60, 64, 38, 26, 65, 84),\\
a_2:=(2, 11, 8, 70, 61, 24, 15, 41, 95, 30, 85, 19, 67, 31, 92, 46, 80, 52, 49, 53, 84, 45, 4, 5, 48,\\
 \ 14, 47, 74, 25, 38, 3, 7, 42, 66, 81, 21, 86, 63, 93, 64, 17, 78, 33, 65, 26, 77, 62, 83, 27, 68,\\ 
 \ 87, 29, 10)(6, 75, 32, 16, 69, 34, 96, 58, 54, 23, 91, 28, 44, 90, 22, 9, 12)\\
 \ (13, 76, 88, 57, 56, 18, 37, 79, 40)(20, 43, 72, 51, 59, 60, 50, 55, 94)(35, 71, 39, 73, 82, 36, 89),\\
a_3:=(3, 18, 55, 81, 41, 90, 11, 48, 92, 23, 62, 85, 30, 94, 53, 67, 26, 31, 22, 75, 40, 58, 49, 84, 63,\\
 \ 52, 95, 9, 35, 72, 69, 46, 78, 28, 73, 7, 4, 17, 56, 82, 42, 89, 12, 47, 91, 24, 61, 86, 29, 93, 54,\\
 \ 68, 25, 32, 21, 76, 39, 57, 50, 83, 64, 51, 96, 10, 36, 71, 70, 45, 77, 27, 74, 8)(5, 37)(6, 38)\\
 \ (13, 20, 87, 43, 16, 34, 60, 65, 79, 14, 19, 88, 44, 15, 33, 59, 66, 80),\\
a_4:=(1, 49, 11, 59)(2, 51, 14, 57, 5, 54, 13, 60, 3, 63, 9, 55, 6, 64, 12, 58, 15, 52, 7, 56, 16, 50,\\
 \ 10, 62, 4, 53, 8, 61)(17, 65, 27, 75)(18, 67, 30, 73, 21, 70, 29, 76, 19, 79, 25, 71, 22, 80, 28,\\
 \ 74, 31, 68, 23, 72, 32, 66, 26, 78, 20, 69, 24, 77)(33, 81, 43, 91)(34, 83, 46, 89, 37, 86, 45,\\
 \ 92, 35, 95, 41, 87, 38, 96, 44, 90, 47, 84, 39, 88, 48, 82, 42, 94, 36, 85, 40, 93),\\
a_5:=(2, 13, 14)(3, 12, 7, 11, 4, 15)(5, 9, 6)(8, 16)(18, 29, 30)(19, 28, 23, 27, 20, 31)(21, 25, 22)\\
 \ (24, 32)(34, 45, 46)(35, 44, 39, 43, 36, 47)(37, 41, 38)(40, 48)(49, 54, 50, 58, 61, 57)\\
 \ (51, 55, 56)(53, 62)(59, 63, 64)(65, 70, 66, 74, 77, 73)(67, 71, 72)(69, 78)(75, 79, 80)\\
 \ (81, 86, 82, 90, 93, 89)(83, 87, 88)(85, 94)(91, 95, 96),\\
a_6:=(1, 11, 9, 3, 13, 12, 6, 8, 14, 16, 10, 7)(2, 15, 5, 4)(17, 27, 25, 19, 29, 28, 22, 24, 30,\\
 \ 32, 26, 23)(18, 31, 21, 20)(33, 43, 41, 35, 45, 44, 38, 40, 46, 48, 42, 39)(34, 47, 37, 36)\\
 \ (49, 63, 58, 51, 50, 59, 54, 52, 61, 64, 53, 56)(55, 62, 60, 57)(65, 79, 74, 67, 66, 75, 70,\\
 \ 68, 77, 80, 69, 72)(71, 78, 76, 73)(81, 95, 90, 83, 82, 91, 86, 84, 93, 96, 85, 88)(87, 94, 92, 89).
$

The orbits $a_1 * U$, $a_2 * U$, $a_3 * U$, $a_4 * U$, $a_5 * U$, $a_6 * U$  have cardinality $288$, $288$, $288$, $48$, $24$ and $24$ 
respectively and the $(96,95)$ - PA $C := a_1 * U \cup a_2 * U \cup a_3 * U \cup a_4 * U \cup a_5 * U \cup a_6 * U$ is $(96, 10)$ - separable. This proves
\begin{theorem}
\label{thm:10MOLS96}
$N(96) \ge 10$.
\end{theorem}

\bigskip 

Now let $n = 108$.
For the finite fields GF$(4)$ and GF(27),  let $z$ and $x$  respectively be primitive elements satisfying   $z^2 = z+1$ and $x^3 = x+2$.  
Let  $F_1$, $F_2$ and $F_3$ be the arrays below.    Also let $V = ((0,v_1), (0,v_2) ... (0,v_{10}))^T$ where 
$(v_1, v_2,  \ldots,   v_{10}) =     (0,    1,    x,    x^2,     x^3=x+2,    x^4 = x^2 + 2x,    
x^5 = 2x^2+x+2,   x^9 = x+1,    x^{25} = 2x^2 + 1,   x^{20} = 2x^2 + x + 1)$.  
Likewise let $W = ((0,w_1), (0,w_2), \ldots,  (0,w_{10}))^T$  
where $(w_1, w_2, \ldots,  w_{10}) = x(v_1, v_2 \ldots,
v_{10}) = (0, x, x^2, x^3, x^4 = x^2+2x, x^5 = 2x^2 + x +2,  x^6 = x^2+x+1,
 x^{10} = x^2+x,  x^{26}=1,    x^{21} = x^2 + 1)$.

The columns of our $(108,10,1)$-DM are then obtained by adding the $9$ vectors 
$aV + bW$     $(a, b \in Z_3)$  to each column of  $F= [F_1 | F_2|F_3]$.

 
\begin{center} 
 $F_1 = \left (
\begin{array}{ccccccccc}
        (0,0)                & (z+1,   0)             & (z,   0)                & (z+1,   0)         \\
        (0,0)                & (0,  x^2)             & (0,   x^2)             & (1,   2x^2)          \\
        (0,0)                & (0,  2x+1)           & (0,   2)                 & (z,   1)             \\
        (0,0)                & (1,  x^2+2x)       & (z+1, x^2+x+1)    & (0,   0)           \\
        (0,0)                & (1,  x^2+x+2)      & (z+1, 1)              & (z+1, 2x^2+x+1)    \\
        (0,0)                & (1,  2x^2+2x+2)  & (z+1, x^2+2)      & (z,   2x^2+1)       \\
        (0,0)                & (z,  2x^2+1)        & (1, x^2+1)          & (1,   x+2)                   \\
        (0,0)                & (z,  2x^2+x+2)    & (1, 2x)                & (z,   x^2)            \\
        (0,0)                & (z+1, x+1)           & (z,  x^2+2x)        & (1,   2x^2+2x)    &      \\
        (0,0)                & (z+1, 2x^2+x)    & (z, 2x^2+2x+1)    & (0,  x^2+x+2)        \\
\end{array}  \right ) , $
\end{center}

\begin{center} 
 $F_2 = \left (
\begin{array}{ccccccccc}
     &(1,   0)          & (0,   0)       &  (z+1,0)          & (z,   0)     \\
     & (1,  2x^2)       & (1,0)          &  (z,2x^2)         & (z,  x^2)    \\
     & (z,   2)         & (z,1)          & (z+1,2x+1)        & (z+1,  0)    \\
     & (z, 2x^2+2)      & (z+1,1)        &   (z,x)           & (0,  x)      \\
     & (0,  2x+1)       & (z,x^2+2x+2)   &   (1,x^2+x+1)     & (z,  2)       \\
     & (z+1,2x^2+x+1)   & (0,x^2+x+1)    &  (0,2)            & (1,  2x^2+1)  \\
     & (z+1,x^2+2x+1)   & (z,x+1)        &  (0,x^2+x)        & (z+1,  2x)    \\
     & (1,  2)          & (z+1,2x^2+x+1) &    (1, 2)         & (0,  x^2)     \\
     & (0,  x+2)        & (z+1, x)       &  (z, 2x^2+2x+2)   & (1,  x^2)     \\
    & (z+1,  0)         & (1,2x^2 + 1)   &  (1, 2x+2)        & (z+1, x^2+2x) \\
\end{array}  \right ) , $
\end{center}

\begin{center} 
 $F_3 = \left (
\begin{array}{ccccccccc}
     &(1,   0)               & (z,   0)             & (1,   0)          & (0,0)        \\
     & (z,   0)             & (z+1,  x^2)          & (z+1, 0)          & (z+1,2x^2)    \\
     & (z+1, 0)               & (1,   x)             & (1,   2)          & (1,x+2)       \\
     & (1,   2x)            & (z,  x)              & (z+1, 2x^2+2x+1)  & (1,2x+1)     \\
     & (0,  2x^2+2x+1)      & (1,  2x^2+2x)        & (z,  x^2+1)       & (z+1,2x^2 +x) \\
     & (z, 2x^2+2x+2)       & (z+1, 0)             & (1,  x^2+1)       & (z,2x^2+2)    \\
     & (1,  1)              & (z,  2x+2)           & (0,  2x+1)        & (z+1,x^2+x+2) \\
     & (z+1,  2x^2+1)       & (z+1,  x^2+1)        & (z,  2x^2+x+2)    & (0,x^2+x+1)   \\
     & (z+1,  x^2+1)        & (1,  2x+1)           & (0, 2x^2+1)       & (z,x+1)       \\
     & (z,  x+2)            & (0,  2x^2+1)         & (z, x)            & (1,2x^2+2x)    \\
 \end{array}  \right ) . $
\end{center}
This proves
\begin{theorem}
	\label{thm:9MOLS108}
	$N(108) \ge 9$.
\end{theorem}
\bigskip
Finaly we come to the correction of an error in \cite{ACD} for the case n=45.
Here we give the corrected $(45,7,1)$ difference matrix. Consider the arrays

\begin{center} 
 $E_1 = \left (
\begin{array}{ccccccccc}
   (0,0,0) \\
   (0,0,0) \\
   (0,0,0) \\
   (0,0,0) \\
   (0,0,0) \\
   (0,0,0) \\
   (0,0,0) \\
\end{array} \right ), \ \ \ E_2 = \left ( \begin{array}{cccccccc}
   (2,2,1) & (3,1,1) & (4,1,2) & (4,0,1) & (0,1,1) & (0,2,1) & (3,2,2) \\ 
   (1,2,1) & (4,2,2) & (1,2,0) & (4,1,0) & (3,1,1) & (3,0,0) & (2,1,2) \\
   (4,1,1) & (2,2,1) & (3,2,0) & (1,2,0) & (2,1,0) & (1,0,0) & (3,2,1) \\
   (0,1,0) & (2,1,1) & (4,0,0) & (0,0,2) & (4,2,2) & (3,2,2) & (1,2,2) \\
   (3,1,2) & (2,1,0) & (0,2,2) & (4,2,1) & (0,2,1) & (2,0,1) & (1,1,2) \\
   (2,1,1) & (1,2,2) & (3,0,1) & (2,0,1) & (1,0,0) & (4,2,1) & (1,1,0) \\
   (0,0,0) & (0,0,0) & (0,0,0) & (0,0,0) & (0,0,0) & (0,0,0) & (0,0,0) \\
\end{array} \right ) . $ \end{center}
  Now let $E_3$ be the array over $Z_5 \times Z_3 \times Z_3$
  obtained by interchanging the following pairs
  of rows in $E_2$: 1 and 6, 2 and 5, 3 and 4,  while leaving row 7  
  unaltered.  Also,  let $V = [(0,0,1),$  $(0 , 2 , 0),$ $(0, 1, 2),$ $(0, 2, 1),$   $(0, 1, 0),$  $(0 , 0, 2 ),$  $(0 , 0 , 0)]^T$. 
The columns of our $(45,7,1)$-DM are then obtained by adding the $3$ vectors $aV$ $(a \in Z_3)$ 
to each column of  $E = [E_1 | E_2 | E_3]$.

\section*{Acknowledgements} 
The authors would like to thank Alice Miller for several useful suggestions and comments,
and for checking a number of results  within this paper.

\section*{ORCID}
R. J. R. Abel:     https://orcid.org/0000-0002-3632-9612\\
I. Janiszczak:      https://orcid.org/0000-0002-4892-1189

\bibliographystyle{plain}


\newpage

	\section*{Appendix:  GAP code}
	\label{appendix:app}
\begin{lstlisting}		

##################################
# HERE YOU NEED TO DEFINE n.
n:=  ;
##################################

a1:=[1..n];
b1:=(1,2);
a2:=[n+1..2*n];
b2:=(n+1,n+2);
t:=Identity(SymmetricGroup([1..2*n]));
for i in [1..n] do
	t:=t*(i,n+i);
od;

# THE GROUPS A,B,B1 AND B2 WILL BE DEFINED
A:=Group(CycleFromList(a1),CycleFromList(a2),b1,b2,t);
B:=Subgroup(A,[CycleFromList(a1),b1,CycleFromList(a2),b2]);
B1:=Subgroup(A,[CycleFromList(a1),b1]);
B2:=Subgroup(A,[CycleFromList(a2),b2]);

###############################################################
# HERE YOU HAVE TO ENTER THE GENERATORS FOR THE ISOMETRY GROUP 
# SEPARATED BY ,
UGenerators:=[];
###############################################################

###############################################################
# HERE ARE THE ORBITREPRESENTATIVES SEPARATED BY ,
OrbReps:=[];
###############################################################

U:=Subgroup(A,UGenerators);

action:=function(x,y) 
#THE ACTION OF ELEMENTS OF A ON ELEMENTS OF B1 
	local x1,x2,xony;
	if x in B then
		x1:=RestrictedPerm(x,a1);
		x2:=RestrictedPerm(x,a2);
		xony:=x1^-1*y*x2^t;
	else
		x1:=RestrictedPerm(x*t,a1);
		x2:=RestrictedPerm(x*t,a2);
		xony:=(x1^-1*y*x2^t)^-1;
	fi;
	return xony;
end;

grouporbit:=function(G,y) 
# RETURNS THE ORBIT GIVEN BY THE ACTION OF G ON y
	local x,orb;
	orb:=[];
	for x in G do
		AddSet(orb,action(x,y));
	od;
	return orb;
end;

dist:=function(x,y) 
# RETURNS THE DISTANCE OF x AND y FOR x,y IN B1
	return(Number(MovedPoints(x^-1 * y)));
end;

mindistset:=function(X) 
# RETURNS THE MINIMUMDISTANCE OF THE SUBSET X OF B1
	local i,j,d,mindist;
	mindist:=n;
	for i in [1..Number(X)] do
		for j in [i+1..Number(X)] do
			d:=dist(X[i],X[j]);
			if d < mindist then
				mindist:=d;
			fi;
		od;
	od;
	return(mindist);
end;		

isdistancebetweensetseq:=function(distance,X,Y) 
# RETURNS true IF FOR ALL x IN X AND ALL y IN Y dist(x,y) EQUALS distance 
# AND false OTHERWISE 
	local i,j;
	for i in [1..Number(X)] do
		for j in [1..Number(Y)] do
			if dist(X[i],Y[j]) = distance then
				continue;
			else
				return(false);
			fi;
		od;
	od;
	return(true);
end;

mkorbits:=function(OrbitReps)
# INPUT IS THE SET OF ORBITREPRESENATIVES
# RETURNS THE CORRESPONDING SET OF ORBITS
	local i, orbs;
	orbs:=[];
	for i in [1..Number(OrbitReps)] do
		orbs[i]:=grouporbit(U,OrbitReps[i]);
	od;
	return(orbs);
end;

mkcode:=function(Orbits)
# INPUT IS A SET OF ORBITS, RETURNS THE UNION OF ALL ORBITS
	local i, Code;
	Code:=[];
	for i in [1..Number(Orbits)] do
		UniteSet(Code,Orbits[i]);
	od;
	return(Code);
end;

nseparable:=function(Code) 
# RETURNS THE Code SEPARATED IN CODES OF LENGTH n WITH MINIMUMDISTANCE n 
# IF THE Code IS (n,Number(Code)/n)-SEPARABLE AND [] OTHERWISE
	local i,j,r, latinsquare, latinsquares, elementshandled,ls1,ls2;
	if Number(Code) mod n > 0 then
		return([]);
	fi;
	r:=Number(Code)/n;
	latinsquares:=[];
	elementshandled:=[];
	for j in [1..Number(Code)] do
		elementshandled[j]:=0;
	od;
	for i in [1..r] do
		latinsquare:=[];
		for j in [1..Number(Code)] do	
			if elementshandled[j] = 1 then
				continue;
			fi;
			if Number(latinsquare) = 0 then
				AddSet(latinsquare,Code[j]);
				elementshandled[j] := 1;
			else
				if dist(latinsquare[1],Code[j]) = n then
					AddSet(latinsquare,Code[j]);
					elementshandled[j] := 1;
				fi;
			fi;
		od;
		if Number(latinsquare) = n then
			if mindistset(latinsquare) = n then
				AddSet(latinsquares,latinsquare);
				continue;
			else
				return([]);
			fi;
		fi;
	od;
	for i in [1..Number(latinsquares)-1] do
		ls1:=latinsquares[i];
		for j in [i+1..Number(latinsquares)] do
			ls2:=latinsquares[j];
			if isdistancebetweensetseq(n-1,ls1,ls2) then
				continue;
			else
				return([]);
			fi;
		od;
	od;	
	return(latinsquares);
end;	

makeMOLS:=function(sepcode)
# GIVEN A (n,m)-SEPARABLE CODE sepcode WHICH IS SEPARATED IN m (n,n) PAs
# THE SET OF m MOLS WILL BE RETURNED
	local i,j,k,sepdeg,MOLS;
	MOLS:=[];
	sepdeg:=Number(sepcode);
	for i in [1..sepdeg] do
		MOLS[i]:=[];
		for j in [1..n] do
			MOLS[i][j]:=[];
		od;
		for j in [1..n] do
			for k in [1..n] do
				MOLS[i][j^sepcode[i][k]][j]:=k;
			od;
		od;
	od;
	return(MOLS);
end;

isls:=function(ls)
# RETURNS true IF ls IS A LATINSQUARE AND false OTHERWISE
	local i,j;
	for i in [1..n] do
		if Number(Set(ls[i])) = n then
			continue;
		else
			return(false);
		fi;
	od;
	return(true);
end;

lsorth:=function(ls1,ls2)
# RETURNS true IF ls1 AND ls2 ARE ORTHOGONAL AND false OTHERWISE
	local i,j,pairset;
	pairset:=[];
	for i in [1..n] do
		for j in [1..n] do
			AddSet(pairset,[ls1[i][j],ls2[i][j]]);
		od;
	od;
	if Number(pairset) = n^2 then
		return(true);
	else
		return(false);
	fi;
end;

orbs:=mkorbits(OrbReps);
# orbs CONTAINS ALL ORBITS FOR THE REPRESENTATIVES IN OrbReps

MOLSCode:=mkcode(orbs);
# MOLSCode IS THE UNION OF ALL ORBITS

SepCode:=nseparable(MOLSCode);
Print(Number(SepCode),"\n");
# PRINTS Number(MOLSCode)/n IF MOLSCode IS SEPARABLE AND 0 ELSE

if Number(SepCode) = Number(MOLSCode)/n then
	Mols:=makeMOLS(SepCode);
fi;
# Mols ARE THE RESULTING MOLS FROM MOLSCODE
	
\end{lstlisting}


\end{document}